\documentclass[11pt]{article}
\usepackage{amsmath}
\usepackage{amssymb}
\usepackage{amsthm}
\usepackage{enumerate}
\usepackage{graphicx,color}
\usepackage{xspace}

\parskip \medskipamount
\newtheorem{theorem}{Theorem}[section]
\newtheorem{lemma}[theorem]{Lemma}
\newtheorem{corollary}[theorem]{Corollary}

\newtheorem{proposition}[theorem]{Proposition}
\newtheorem{remark}[theorem]{Remark}
\newtheorem{definition}[theorem]{Definition}

\def\calC{{\mathcal C}}
\def\calL{{\mathcal L}}

\newcommand{\pr}[1]{\operatorname{\mathbf{P}}\left(#1\right)}

\newcommand{\E}[1]{\operatorname{\mathbf{E}}\left[#1\right]}

\newcommand{\Econd}[2]{\operatorname{\mathbf{E}}\left[#1\;\middle\vert\;#2\right]}

\newcommand{\critical}{\mathrm{c}}

\newcommand{\tv}[2]{\|#1-#2\|_\mathrm{TV}}
\newcommand{\tmix}{T_\mathrm{mix}}

\newcommand{\fresh}{\star}
\newcommand{\regeneration}{\Omega_\mathrm{REG}}

\newcounter{notecounter}

\newcommand{\IGNORE}[1]{}

\newcommand{\ndtorus}{\mathbb{T}^{d,n}}
\newcommand{\bfZ}{\mathbb{Z}^d}
\newcommand{\dntorus}{\mathbb{T}^{d,n}}
\newcommand{\twontorus}{\mathbb{T}^{2,n}}
\newcommand{\onentorus}{\mathbb{T}^{1,n}}
\newcommand{\dist}{\mathrm{dist}}

\newcommand{\tmixRW}{T^{\mathrm{RW}}_{\mathrm{mix}}}
\newcommand{\Pruu}[1]{\operatorname{\mathbf{P}}\left[#1\right]}
\newcommand{\bfP}{\mathbf{P}}

\newtheorem{remarks}[theorem]{Remarks}
\def\calN{{\mathcal N}}
\begin{document}

\title{An elementary proof of Marton and Shields' obstruction to finitary coding
}

\author{
        Jeffrey E. Steif\thanks{Chalmers University of Technology
and Gothenburg University, Gothenburg, Sweden\ \ Email:
        \hbox{steif@chalmers.se}}
}

\maketitle
\thispagestyle{empty}



  \medskip\noindent


Recall that if $X=\{X_i\}_{i\in Z^d}$ and $Y=\{Y_i\}_{i\in Z^d}$ are stationary fields taking values in
$U$ and $V$ respectively, we say that $Y$ is a factor of $X$ if there is a mapping $f$ from
$U^{Z^d}$ to $V^{Z^d}$ which is equivariant (commutes with shifts in $Z^d$) such that the image of
$X$ under $f$ is $Y$. $f$ is then called a factor map. Furthermore, we say that $Y$ is a finitary factor of
$X$ if $f$ can be taken to have the property that a.s.\ $Y_0$ is determined by the $X$ process in a sufficiently large 
(random) box around 0. $f$ is then called a finitary map.

The main result concerning obstructions to finitary coding from an i.i.d.\ process
is due to K. Marton and P. Shields in their article
The positive-divergence and blowing up properties, appearing in Israel Journal of Mathematics, volume 86, pages
331-348, 1994.

They showed that (1) a finitary factor of an i.i.d. process satisfies the {\it blowing up property} which I won't bother to define and which we don't need, and (2) the blowing up property implies that one has 
{\it usual large deviation behavior} in the ergodic theorem
meaning that for fixed $\epsilon>0$, the probability that the average of the process
between time 1 and $n$ differs from the mean by more than $\epsilon>0$ decays exponentially in $n$.

Lots of processes don't have  "usual large deviation behavior" and therefore cannot 
be a finitary factor of an i.i.d. process by the Marton-Shields Theorem. (For example, the plus state for the Ising model
in the coexistence regime or the range of a random walk run for time $n$.)
This is the only way I know how to show something which is a factor of i.i.d.'s is not however a finitary factor
and then one does not need to deal with the so-called blowing up property which is nice.
The above is essentially true. Of course if you happen to have a process which does satisfy 
"usual large deviation behavior" but not the blowing up property, then the above would not be true but I am
not aware of any (natural) examples like that.

When I read the Marton-Shields paper (maybe 25 years ago or more), 
I wondered how hard it would be to prove the {\bf weaker} statement
(which is of course an immediate corollary of the Marton-Shields theorem)
that  a finitary factor of an i.i.d. process has "usual large deviation behavior". I managed to do a simple proof
of this {\bf weaker} statement and I have meant to put it on my homepage for years
(since obviously this is not publishable) and finally I am doing this now.

\newpage

Theorem:  Let  $X=\{X_n\}_{n\in Z^d}$ be an i.i.d.\ process and 
let $\{Y_n\}_{n\in Z^d}$  be a process taking values in a bounded interval which is
a finitary factor of $X$. Then $\{Y_n\}_{n\in Z^d}$  has usual large deviation behavior.

Proof:
 Without loss of generality, we take $d=1$ (the same argument works in higher dimensions),
the $Y$ process takes values in $[-1,1]$ and  $E(Y_0)=0$. 
Let  $S_n:=\sum_{i=1}^n Y_i$. We want to show that for every $\epsilon >0$,
$P(\frac{|S_n|}{n}\ge \epsilon)$ decays exponentially in $n$.

Fix $\epsilon >0$. Clearly one 
can choose $\delta >0$ such that for any finite sequence $Y_1,\ldots,Y_n$ taking
values in $[-1,1]$, if  $\tilde{Y_1},\ldots,\tilde{Y_n}$ differs from 
 $Y_1,\ldots,Y_n$ in at most $\delta n$ locations and takes the value zero at these locations where  they differ, then
 $$
 | (Y_1+\dots+Y_n)/n-  (\tilde{Y_1}+\dots+\tilde{Y_n})/n| \le \frac{\epsilon}{4}.
$$

Now let $N_i$ be the diameter of the smallest box around $i$ in the $X$ process
that determines $Y_i$. So $N_i$ is a random variable
and is finite a.s. Now, choose $N$ so large that
$$
P(N_0\ge N)<  \frac{\delta}{2}
$$
and  (using that $E(Y_0)=0$ and dominated convergence)
$$
|E[Y_0 I_{\{ N_0\le N\}}]| \le \frac{\epsilon}{4}
$$

Let $Z_i:=I_{\{N_i > N\}}$  and let
$\tilde{Y_i}:=Y_i I_{\{N_i\le N\}}$. Note that the processes $\{Z\}$ and 
$\{\tilde{Y}\}$ are both the image of a block (finite range) map of an i.i.d.\ process
 and hence satisfy standard large deviation behavior
(using for example Azuma's Lemma).

By definition of $\delta$, we have that if $\sum_{i=1}^n Z_i\le \delta n$, then
 $$
 | (Y_1+\dots +Y_n)/n-  (\tilde{Y_1}+\dots+\tilde{Y_n})/n| \le \frac{\epsilon}{4}.
$$

It follows that 
$$
P(\frac{|S_n|}{n}> \epsilon)\le 
P(\sum_{i=1}^n Z_i\ge \delta n)+
P(|\tilde{Y_1}+\dots+\tilde{Y_n})/n| \ge \frac{3\epsilon}{4}).
$$

By the first condition on how we chose $N$, we have that the $Z_i$'s have mean at most $\delta/2$
and since that process is the image of a block map, the first summand goes to 0 exponentially.
By the second condition on how we chose $N$, we have that $|E(\tilde{Y_0})|\le \frac{\epsilon}{4}$.
Since that process is also the image of a block map, the second summand also goes to 0 exponentially.\\
QED

\bigskip
{\it Remark}  \, Often the definition of finitary is taken to mean that the factor map is continuous a.s.
This is equivalent to the definition given above in terms of a "stopping box" when the sets $U$ and $V$ are
countable with the discrete topology. However, when one goes beyond these cases, then these two
definitions are no longer equivalent and then one has to be careful about what one means by a 
 finitary map.

\bigskip
{\it Acknowledgements} Thanks to Yinon Spinka for various comments and in particular for pointing out to me
the contents of  the final remark in the paper and the suggestion of Azuma's Lemma.

\end{document}